\documentstyle[12pt]{article}
\def\a{\alpha} 
\def\b{\beta}

\def\D{\triangle}
\def\t{\tau}
\def\d{\delta}  
\def\th{\theta} 
\def\l{\lambda}
\def\L{\Lambda}

\def\s{\sigma}
\def\sou{\overline}
\def\so{\underline} 
\def\ch{\hat} 

\def\f{\rightarrow}

\def\p{\succ}

\def\N{\ifmmode{\rm I\mkern-3.1mu
N\mkern0.5mu}\else{\rm I\kern-.16em
N\hskip0.5pt\ }\fi\relax}

\begin{document} 

\Large\bf 
$\so{S}$-STORAGE OPERATORS\\
 
\normalsize \bf 
Karim NOUR \footnote{We wish to thank Ren\'e David for helpful discussions.}\\
\rm
LAMA - Equipe de Logique, Universit\'e de Savoie -
73376 Le Bourget du Lac cedex \footnote{e-mail nour@univ-savoie.fr} \\

\bf Abstract \rm 
In 1990, J.L. Krivine introduced the notion of storage operator to simulate, for Church integers,
the ``call by value'' in a context of a ``call by name'' strategy. In this present paper, we
define, for every $\l$-term $\so{S}$ which realizes the successor function on Church integers, the
notion of $\so{S}$-storage operator. We prove that every storage operator is a
$\so{S}$-storage operator. But the converse is not always true.  \\

\bf Mathematics Subject Classification \rm: 03B40, 68Q60 \\
\bf Keywords \rm: Church integer ; Storage operator ; Call by value ; Call by name ; Head
reduction ; Solvable ; Successor ; $\so{S}$-storage operator.

\section{Definitions and notations}

\begin{itemize}
\item We denote by $\L$ the set of $\l$-terms modulo $\a$-equivalence, and by $\cal V$ the set of
$\l$-variables. 
\item Let $t,u,u_1,...,u_n$ be $\l$-terms, the application of $t$ to $u$ is denoted
by $(t)u$. In the same way we write $(t)u_1...u_n$ instead of $(...((t)u_1)...)u_n$. 
\item The sequence of $\l$-terms $u_1,...,u_n$ is denoted $\sou{u}$. 
\item If $\sou{u} = u_1,...,u_n$, we denote by $(t)\sou{u}$ the $\l$-term $(t)u_1...u_n$.   
\item The $\b$-equivalence relation is denoted by $u \simeq_{\b} v$.  
\item The notation $\s(t)$ represents the result of the simultaneous
substitution $\s$ to the free variables and the constants of $t$ after a suitable renaming of the
bounded variables of $t$.   
\item Let us recall that a $\l$-term $t$ either has a head redex [i.e. $t=\l x_1 ...\l
x_n (\l x u) v \sou{w}$, the head redex being $(\l x u) v$], or is in head normal form
[i.e. $t=\l x_1 ...\l x_n (x) \sou{w}$].  
\item The notation $u \p v$ means that $v$ is obtained from $u$ by some head reductions.
\item If $u \p v$, we denote by $h(u,v)$ the length of the head reduction between $u$ and $v$.  
\item A $\l$-term $t$ is said solvable iff the head reduction of $t$ terminates.
\end{itemize}

The following results are well known (see [3]):
\begin{itemize}
\item[] --{\it If $u$ is $\b$-equivalent to a solvable $\l$-term, then $t$ is solvable.}
\item[] --{\it If $u \p v$, then, for any substitution $\s$, $\s(u) \p \s(v)$, and
$h(\s(u),\s(v))$=h(u,v).\\ In particular, if for some  substitution $\s$, $\s(t)$ is solvable,
then $t$ is solvable.} 
\end{itemize}

\begin{itemize}
\item We define $(u)^n v$ by induction : $(u)^0 v = v$ and $(u)^{n+1} v = (u)(u)^n v$.
\item For each integer $n$, we define the Church integer $\so{n} = \l f\l x(f)^n x$.
\item A closed $\l$-term $\so{S}$ is called successor iff, for every $k \geq 0$, $(\so{S})\so{k}
\simeq_{\b} \so{k+1}$.   
\end{itemize}

{\bf Examples} Let $\so{S_1}=\l n \l f \l x (f)((n)f)x$ and $\so{S_2}=\l n \l f \l x ((n)f)(f)x$.
\\It is easy to check that $\so{S_1}$ and $\so{S_2}$ are successors. $\Box$

\section{Introduction}

In $\l$-calculus the left reduction strategy (iteration of the head reduction) has much advantages
: it always terminates when applied to a normalizable $\l$-term and it seems more economic since
we compute a $\l$-term only when we need it. But the major drawback of this strategy is that a
function must compute its argument every time it uses it. In 1990 J-L. Krivine introduced the
notion of storage operators in order to avoid this problem and to simulate call-by-value when
necessary.\\

Let $F$ be a $\l$-term (a function), and $\so{n}$ a Church integer. During the computation, by left
reduction, of $(F)\th_n$ (where $\th_n \simeq\sb{\b} \so{n}$), $\th_n$ may be computed several times
(as many times as $F$ uses it). We would like to transform $(F)\th_n$ to $(F)\t_n$ where $\t_n$ is a
fixed closed $\l$-term $\b$-equivalent to $\so{n}$. We also want this transformation depends only on
$\th_n$ (and not $F$). \\

Therefore the definition : A closed $\l$-term $T$ is called storage operator if and
only if for every $n \geq 0$, there is a closed $\l$-term $\t_n \simeq\sb{\b} \so{n}$ such that for
every $\th_n \simeq\sb{\b} \so{n}$, $(T)\th_n f \p (f)\t_{n}$ (where $f$ is a new variable). \\ 

Let's analyse the head reduction $(T)\th_n f \p (f)\t_n$, by replacing each $\l$-term which comes
from $\th_n$ by a new variable.\\

If  $\th_n \simeq\sb{\b} \so{n}$,  then $\th_n \p \l g \l x(g)t_{n-1}$,  $t_{n-k} \p (g)t_{n-k-1}$ 
$1 \leq k \leq n-1$,  $t_0 \p x$,  and  $t_k \simeq\sb{\b} (g)^k x$ $0 \leq k \leq n-1$.  \\

Let $x_n$ be a new variable ($x_n$ represents $\th_n$). $(T)x_nf$ is solvable, and its head normal
form does not begin by $\l$, therefore it is a variable applied to some arguments. The free
variables of $(T)x_n f$ are $x_n$ and $f$, we then have two possibilities for its head normal form :
$(f)\d_n$ (in this case we stop) or $(x_n)a_1...a_m$. \\

Assume we obtain $(x_n)a_1...a_m$. The variable $x_n$ represents $\th_n$, and  $\th_n \p \l g \l
x(g)t_{n-1}$, therefore $(\th_n)a_1...a_m$ and $((a_1)t_{n-1}[a_1/x,a_2/g])a_3...a_m$ have the
same head normal form. The $\l$-term $t_{n-1}[a_1/g,a_2/x]$ comes from $\th_n$. Let
$x_{n-1,a_1,a_2}$ be a new variable $(x_{n-1,a_1,a_2}$ represents $t_{n-1}[a_1/g,a_2/x])$. The
$\l$-term $((a_1)x_{n-1,a_1,a_2})a_3...a_m$ is solvable, and its head normal form does not begin
by $\l$, therefore it is a variable applied to some arguments. The free variables of
$((a_1)x_{n-1,a_1,a_2})a_3...a_m$ are among $x_{n-1,a_1,a_2}$, $x_n$, and $f$, we then have three
possibilities for its head normal form :  $(f)\d_n$ (in this case we stop) or $(x_n)b_1...b_r$ or
$(x_{n-1,a_1,a_2})b_1...b_r$.\\

Assume we obtain $(x_{n-1,a_1,a_2})b_1...b_r$. The variable $(x_{n-1,a_1,a_2}$ represents
$t_{n-1}[a_1/g,a_2/x]$, and $t_{n-1} \p (g)t_{n-2}$, therefore $(t_{n-1}[a_1/g,a_2/x])b_1...b_r$
and $((a_1)t_{n-2}[a_1/g,a_2/x])b_1...b_r$ have the same head normal form. The $\l$-term
$t_{n-2}[a_1/g,a_2/x]$ comes from $\th_n$. Let $x_{n-2,a_1,a_2}$ be a new variable
$(x_{n-2,a_1,a_2}$ represents $t_{n-2}[a_1/g,a_2/x])$. The $\l$-term
$((a_1)x_{n-2,a_1,a_2})b_1...b_r$ is solvable, and its head normal form does not begin by $\l$,
therefore it is a variable applied to arguments. The free variables of
$((a_1)x_{n-2,a_1,a_2})b_1...b_r$ are among $x_{n-2,a_1,a_2}$, $x_{n-1,a_1,a_2}$, $v_n$, and $f$,
therefore we have four possibilities for its head normal form : $(f)\d_n$ (in this case we stop)
or $(x_n)c_1...c_s$ or $(x_{n-1,a_1,a_2})c_1...c_s$ or $(x_{n-2,a_1,a_2})c_1...c_s$ ... and so
on...\\ 

Assume we obtain $(x_{0,d_1,d_2})e_1...e_k$ during the construction. The variable
$x_{0,d_1,d_2}$ represents $t_0[d_1/g,d_2/x]$, and $t_0 \p x$, therefore
$(t_0[d_1/g,d_2/x])e_1...e_k$ and $(d_2)e_1...e_k$ have the same head normal form ; we then
follow the construction with the $\l$-term $(d_2)e_1...e_k$. The $\l$-term $(T)\th_n f$ is
solvable, and has $(f)\t_n$ as head normal form, so this construction always stops on $(f)\d_n$.
We can prove by a simple argument that $\d_n \simeq\sb{\b} \so{n}$.  \\

According to the previous construction, the reduction $(T)\th_nf \p (f)\t_n$ can
be divided into two parts :  
\begin{itemize}
\item[] - A reduction that does not depend on $n$ : 
\begin{eqnarray} 
(T)x_nf &\p &(x_n)a_1...a_m \nonumber \\
((a_1)x_{n-1,a_1,a_2})a_3...a_m &\p &(x_{n-1,a_1,a_2})b_1...b_r \nonumber \\
((a_1)x_{n-2,a_1,a_2})b_1...b_r &\p &(x_{n-2,a_1,a_2})b_1...b_r \nonumber \\
&. \nonumber \\
&. \nonumber \\
&. \nonumber
\end{eqnarray} 
\item[] - A transformation that depends on $n$ (and not on $\th_n$) :  
\begin{eqnarray}
(x_n)a_1...a_m &\leadsto &((a_1)x_{n-1,a_1,a_2})a_3...a_m  \nonumber \\
(x_{n-1,a_1,a_2})b_1...b_r &\leadsto &((a_1)x_{n-2,a_1,a_2})c_1...c_s \nonumber \\
&.\nonumber \\
&.\nonumber \\
&. \nonumber \\  
(x_{0,d_1,d_2})e_1...e_k &\leadsto &(d_1)e_1..e_k \nonumber 
\end{eqnarray}
\end{itemize}

We add new constants $x_i$ and $x_{i,a,b,\sou{c}}$ in $\l$-calculus, and we consider the following set
of head reduction rules : 
\begin{eqnarray}
(\l x u)v \sou{w} &\p &(u[v/x])\sou{w} \nonumber \\
(x_{i+1})ab\sou{c} &\p &((a)x_{i,a,b,\sou{c}})\sou{c} \nonumber \\
(x_{0})ab\sou{c} &\p &(b)\sou{c} \nonumber \\
(x_{i+1,a,b,\sou{c}})\sou{w} &\p &((a)x_{i,a,b,\sou{w}})\sou{w} \nonumber \\
(x_{0,a,b,\sou{c}})\sou{w} &\p &(b)\sou{w} \nonumber 
\end{eqnarray}
We write $t \p_x t'$ if $t'$ is obtained from $t$ by applying these rules finitely many times.\\

With this formalisme we have the following result (see [1] and [4]):\\
{\it A closed $\l$-term $T$ is a storage operator iff for every $n \geq 0$, $(T)x_n f \p_x (f)\t_n$
and where $\t_n$ is a closed $\l$-term $\b$-equivalent to $\so{n}$.}\\

The constants $x_i$ and $x_{i,a,b,\sou{c}}$ represent intuitively the $\l$-terms which come
from a non calculated Church integer. The uniform shape of Church integers allows to
describe the behaviour of these constants when they are in the head position. However, another
method to describe a Church integer is simply to say that it is zero or a successor.\\

Formally, we add new constants $X_i$ et $X_{i,a,b,\sou{c}}$ in $\l$-calculus, and we consider, for
every successor $\so{S}$, the following set of head reduction rules :
\begin{eqnarray}
(\l x u)v \sou{w} &\p &(u[v/x])\sou{w} \nonumber \\
(X_{i+1})ab\sou{c} &\p &((\so{S})X_{i,a,b,\sou{c}})ab\sou{c} \nonumber \\
(X_{0})ab\sou{c} &\p &(\so{0})ab\sou{c} \nonumber \\
(X_{i+1,a,b,\sou{c}})uv\sou{w} &\p &((\so{S})X_{i,u,v,\sou{w}})uv\sou{w} \nonumber \\
(X_{0,a,b,\sou{c}})uv\sou{w} &\p &(\so{0})uv\sou{w} \nonumber
\end{eqnarray} 
We  write $t \p_X t'$ if $t'$ is obtained from $t$ by applying these rules finitely many times.\\

A $\so{S}$-storage operator is defined as follows :\\
A closed $\l$-term $T$ is a $\so{S}$-storage operator iff for every $n \geq 0$, $(T)X_n f \p_X
(f)\t_n$ where $\t_n$ is a closed $\l$-term $\b$-equivalent to $\so{n}$.\\

This paper studies the link  betwen the storage operators and the $\so{S}$-storage operators. 
We prove that every storage operator is a $\so{S}$-storage operator. But the converse is not
always true. 

\section{Storage operators and $\so{S}$-storage operators}

{\bf Definition} Let $T$ be a closed $\l$-term. We say that $T$ is a
storage operator iff for every $n \geq 0$, there is a closed \footnote{In his definition
of storage operator, J.L. Krivine autorizes the $\t_n$ to contain free variables which are
replaced by terms depend of $\th_n$. The results of this paper remain valid with this
definition but the proofs will be too technical.} $\l$-term $\t_n \simeq\sb{\b}
\so{n}$, such that for every  $\th_n \simeq\sb{\b} \so{n}$, $(T)\th_n f \p (f)\t_n$. \\

{\bf Remark} Let $F$ be any $\l$-term (for a function), and $\th_n$ a $\l$-term $\b$-equivalent to
$\so{n}$. During the computation of $(F)\th_n$, $\th_n$ may be computed each time it comes in head
position. Instead of computing $(F)\th_n$, let us look at the head reduction of $(T)\th_n F$. Since
it is $\{(T)\th_n f\}[F/f]$, we shall first reduce $(T)\th_n f$ to its head normal form,
which is $(f)\t_n$, and then compute $(F)\t_n$. The computation has been decomposed into two
parts, the first being independent of $F$. This first part is essentially a computation of $\th_n$,
the result being $\t_n$, which is a kind of normal form of $\th_n$. So, in the computation of
$(T)\th_n F$, $\th_n$ is computed first, and the result is given to $F$ as an argument, $T$ has
stored the result, before giving it, as many times as needed, to any function. $\Box$ \\

{\bf Examples} Let $\so{S}$ be a  successor. If we take :\\ 
$T_1 = \l n((n)G)\d$ where
$G = \l x\l y(x)\l z(y)(\so{S})z$ and $\d = \l f(f)\so{0}$ \\
$T_2 = \l n\l f(((n)F)f)\so{0}$
where $F = \l x\l y(x)(\so{S})y$, \\
then it is easy to check that (see [1] and [3]): \\
for every $\th_n \simeq\sb{\b} \so{n}$, $(T_i)\th_n f \p (f)(\so{S})^n \so{0}$ ($i=1$ or $2$). \\ 
Therefore $T_1$ and $T_2$ are storage operators. $\Box$ \\

 Let $\{x_i\}_{i\geq 0}$ be a set of different constants.
We define a set of terms (denoted by $\L_x$) in the following way :
\begin{itemize}  
\item[] - If $x\in {\cal V}\bigcup\{x_i\}_{i\geq 0}$, then $x \in \L_x$ ;
\item[] - If $x \in {\cal V}$, and $u \in \L_x$, then $\l xu \in \L_x$ ;
\item[] - If $u \in \L_x$, and $v \in \L_x$, then $(u)v \in \L_x$ ;
\item[] - If $n \in \N$, and $a,b,\sou{c} \in \L_x$, then $x_{n,a,b,\sou{c}} \in \L_x$. \\
\end{itemize}
$x_{n,a,b,\sou{c}}$ is considered as a constant which does not appear in $a,b,\sou{c}$.\\
The terms of the set $\L_x$ are called $\l x$-terms.\\

We have the following result (see [1] and [4]) :\\
{\it A closed $\l$-term $T$ is a storage operator
iff for every $n \geq 0$, there is a finite sequence of head reduction $\{ U_i \p V_i \}_{1\leq i\leq
r}$ such that : \begin{itemize}
\item[] 1) $U_i$ and $V_i$ are $\l x$-terms ; 
\item[] 2) $U_1 = (T)x_n f$ and $V_r = (f)\t_n$  where $\t_n$ is closed $\l$-term
$\b$-equivalent to $\so{n}$ ;  
\item[] 3) $V_i = (x_n) a b \sou{c}$ or $V_i = (x_{l,a,b,\sou{c}}) \sou{d}$  $0\leq l\leq n-1$;
\item[] 4) If $V_i = (x_n)a b \sou{c}$, then $U_{i+1} = (b)\sou{c}$ if $n=0$ and
$U_{i+1} = ((a)x_{n-1,a,b,\sou{c}})\sou{c}$ if $n \neq 0$ ;
\item[] 5) If $V_i = (x_{l,a,b,\sou{c}})\sou{d}$ $0\leq l\leq n-1$, then $U_{i+1} = (b)\sou{d}$ if
$l=0$ and $U_{i+1} = ((a)x_{l-1,a,b,\sou{d}})\sou{d}$ if $l \neq 0$.\\
\end{itemize} }

{\bf Definitions} \\
{\bf 1)} Let $\{X_i\}_{i\geq 0}$ be a set of different constants.
We define a set of terms (denoted by $\L_X$) in the following way :
\begin{itemize}  
\item[] - If $x\in {\cal V}\bigcup\{X_i\}_{i\geq 0}$, then $x \in \L_X$ ;
\item[] - If $x \in {\cal V}$, and $u \in \L_X$, then $\l xu \in \L_X$ ;
\item[] - If $u \in \L_X$, and $v \in \L_X$, then $(u)v \in \L_X$ ;
\item[] - If $n \in \N$, and $a,b,\sou{c} \in \L_X$, then $X_{n,a,b,\sou{c}} \in \L_X$. \\
\end{itemize}
$X_{n,a,b,\sou{c}}$ is considered as a constant which does not appear in $a,b,\sou{c}$. \\
The terms of the set $\L_X$ are called $\l X$-terms.\\

{\bf 2)} Let $\so{S}$ be a successor. A closed $\l$-term $T$ is called a
$\so{S}$-storage operator iff for every $n \geq 0$, there is a finite sequence of head reduction $\{
U_i \p V_i \}_{1\leq i\leq r}$ such that :
\begin{itemize}
\item[] 1) $U_i$ and $V_i$ are $\l X$-terms ;
\item[] 2) $U_1 = (T)X_n f$ and $V_r = (f)\t_n$ where $\t_n$ is closed $\l$-term
$\b$-equivalent to $\so{n}$ ;  
\item[] 3) $V_i = (X_n) a b \sou{c}$ or $V_i = (X_{l,a,b,\sou{c}}) u v \sou{w}$ 
$0\leq l\leq n-1$;
\item[] 4) If $V_i = (X_n)a b \sou{c}$, then $U_{i+1} = (\so{0})a b \sou{c}$ if $n=0$ and
$U_{i+1} = ((\so{S})X_{n-1,a,b,\sou{c}})a b\sou{c}$ if $n \neq 0$ ;
\item[] 5) If $V_i = (X_{l,a,b,\sou{c}})u v\sou{w}$ $0\leq l\leq n-1$, then $U_{i+1} = (\so{0})u v
\sou{w}$ if $l=0$ and $U_{i+1} = ((\so{S})X_{l-1,u,v,\sou{w}})u v\sou{w}$ if $l \neq 0$.
\end{itemize}

{\bf Examples} It is easy to check that, for $1 \leq i,j \leq 2$, the
above operator $T_i$ is an $\so{S_j}$-storage operator. We check here
(for example) that $T_2$ is an  $\so{S_2}$-storage operator:

Let $n$ be an integer.

If $n = 0$, then we check that $(T_2)X_nf \p (X_n) F ~ f ~ \so{0}$ and
$(\so{0}) ~ F ~ f ~ \so{0} \p (f)\so{0}$.

If $n \neq 0$, then we check that:
\begin{eqnarray}
(T)X_nf &\p & (X_n) ~ F ~ f ~ \so{0} \nonumber\\
((\so{S_2})X_{n-1,F,f,\so{0}}) ~ F ~ f ~ \so{0}
&\p &(X_{n-1,F,f,\so{0}}) ~ F ~ (F)f ~ \so{0} \nonumber \\
\, &. &\, \nonumber\\
\, &. &\, \nonumber\\
\, &. &\, \nonumber\\
((\so{S_2})X_{0,F,(F)^{n-1}f,\so{0}}) ~ F ~ (F)^{n-1}f ~ \so{0}
&\p &(X_{0,F,(F)^{n-1}f,\so{0}}) ~ F ~ (F)^nf ~ \so{0} \nonumber \\
(\so{0}) ~ F ~ (F)^nf ~ \so{0} &\p &(F)^nf ~ \so{0} \nonumber 
\end{eqnarray}

We prove (by induction on $k$) that, for every $\l$-term $u$, and for every $0 \leq k \leq n$, we
have $(F)^kf ~ u  \p (f)(\so{S_2})^k u$.  
\begin{itemize}
\item[] - For $k=0$, it is true. 
\item[] - Assume that is true for $k$, and prove it for $k+1$.\\ 
$(F)^{k+1}f ~ u = (F) (F)^k f ~ u \p (F)^k f ~ (\so{S_2}) u$. By induction hypothesis we have that for
every $\l$-term $v$, $(F)^k f ~ v \p (f)(\so{S_2})^k v$, then $(F)^{k+1}f ~ u \p 
(f)(\so{S_2})(\so{S_2})^{k} u = (f)(\so{S_2})^{k+1} u$.
\end{itemize}
In particular, for $u=\so{0}$ and $k=n$, we have $(F)^nf ~ \so{0} \p (f)(\so{S_2})^n \so{0}$.\\
Therefore $T_2$ is a $\so{S_2}$-storage operator. $\Box$ \\

A question arizes : {\bf Is there a link between the storage operators and the
$\so{S}$-storage operators ?} 

\section{Link between the storage operators and the $\so{S}$-storage operators}

{\bf Theorem 1} {\it If $T$ is a storage operator, then, for every successor $\so{S}$, $T$ is a
$\so{S}$-storage operator.}\\

{\bf Proof} Let $\so{S}$ be a successor and $T$ a storage operator.\\
Then for every $n \geq 0$, there is a closed $\l$-term $\t_n \simeq_{\b}\so{n}$ such that for every
$\th_n \simeq_{\b}\so{n}$, $(T)\th_n f \p (f)\t_n$. In particular $((T)(\so{S})^n\so{0}) f \p
(f)\t_n$.\\ Let $\s : \L_X \f \L$ the simultaneous substitution defined by :
\begin{center}
$\s(X_n)=(\so{S})^n\so{0}$\\
for every $0\leq k\leq n-1$, $\s(X_{k,a,b,\sou{c}})=(\so{S})^k\so{0}$\\
$\s(x)=x$ if $x \neq X_n , X_{k,a,b,\sou{c}}$
\end{center}
 For every $n \geq 0$, we construct a set of head equation  $\{ U_i \p V_i \}_{1\leq i\leq r}$ such
that : \begin{itemize}
\item[] 1) $U_i$ and $V_i$ are $\l X$-terms ;
\item[] 2) $V_r=(f)\d_n$ ; 
\item[] 3) for every $1\leq i\leq r-1$, $V_i=(X_n)ab\sou{c}$ or
$V_i=(X_{l,a,b,\sou{c}})uv\sou{w}$ ; \item[] 4) $\s(V_i)$ is solvable.
\end{itemize}

Let $U_1 = (T)X_n f$. We have $\s(U_1) = ((T)(\so{S})^n\so{0})f$ is solvable, then $U_1$
is solvable and $U_1 \p V_1$ where $V_1 = (f)\d_n$ or $V_1 = (X_n)ab\sou{c}$. It is clear that
$\s(V_1)$ is solvable. \\
Assume that we have the head reduction $U_k \p V_k$ and $V_k \neq (f)\d_n$. 
\begin{itemize} 
\item[] 
- If $V_k = (X_n)a b \sou{c}$, then, by induction hypothesis, 
$\s(V_k)=((\so{S})^n)\so{0})\s(a) \s(b) \sou{\s(c)}$ is solvable.  
\begin{itemize}
\item[]
- If $n=0$, let  $U_{k+1} = (\so{0})ab\sou{c}$. Then $\s(U_{k+1})=(\so{0})\s(b)\s(b)\sou{\s(c)}$ is
solvable. \item[] 
- If $n \neq 0$, let $U_{k+1} = ((\so{S})X_{n-1,a,b,\sou{c}})ab\sou{c}$. Then
$\s(U_{k+1})=$ \\$((\so{S})(\so{S})^{n-1})\so{0})\s(a)\s(b)\sou{\s(c)}=\s(V_k)$ is solvable.  
\end{itemize}  
\end{itemize}
\begin{itemize}
\item[] 
- If $V_k = (X_{l,a,b,\sou{c}})uv \sou{w}$, then, by induction hypothesis,
$\s(V_k)=((\so{S})^l)\so{0})\s(u) \s(v) \sou{\s(w)}$ is solvable.  
\begin{itemize}
\item[] 
- If $l=0$, let $U_{k+1}=(\so{0})uv\sou{w}$. Then $\s(U_{k+1})=(\so{0})\s(u)\s(v)\sou{\s(w)}$ is
solvable. \item[] 
- If $l \neq 0$,  let $U_{k+1} = ((\so{S})X_{l-1,u,v,\sou{w}})uv\sou{w}$. 
Then $\s(U_{k+1})=$\\$((\so{S})(\so{S})^{l-1})\so{0})\s(u)\s(v)\sou{\s(w)}=\s(V_k)$ is solvable.  
\end{itemize}
\end{itemize}
Therefore $U_{k+1}$ is solvable and $U_{k+1} \p V_{k+1}$ where $V_{k+1} = (f)\d_n$ or 
$V_{k+1} = (X_n)a'b'\sou{c'}$ or $V_{k+1} = (X_{r,a',b',\sou{c'}})a''b''\sou{c''}$. Since
$\s(U_{k+1})$ is solvable, then $\s(V_{k+1})$ is also solvable. \\   

This constraction always terminates (i.e there is a $r\geq 0$ such that $V_r = (f)\d_n$). Indeed, if
not, we check easily that the $\l$-term $((T)(\so{S})^n\so{0})f$ is not solvable. \\ 

Let $y$ be a variable, $\so{\ch{S}}=(\l x \so{S})y$, and $\so{\ch{0}}=(\l x \so{0})y$.\\ 
Let $\ch{\s} : \L_X \f \L$ the simultaneous substitution defined by :
\begin{center}
$\ch{\s}(X_n)=(\so{\ch{S}})^n\so{\ch{0}}$\\
for every $0\leq k\leq n-1$, $\ch{\s}(X_{k,a,b,\sou{c}})=(\so{\ch{S}})^k\so{\ch{0}}$\\
$\ch{\s}(x)=x$ if $x \neq X_n,X_{k,a,b,\sou{c}}$
\end{center}
Since $(\so{\ch{S}})t \p (\so{S})t$ and $\so{\ch{0}} \p \so{0}$, we check easily that 
$((T)(\so{\ch{S}})^n\so{\ch{0}})f \p (f) \ch{\s}(\d_n)$. But $(\so{\ch{S}})^n\so{\ch{0}} \simeq_{\b}
\so{n}$, then  $((T)(\so{\ch{S}})^n\so{\ch{0}})f \p (f) \t_n$. Therefore $\ch{\s}(\d_n)=\t_n$.
Since $\t_n$ is closed, then  $\d_n$ is also closed and $\d_n=\t_n \simeq_{\b} \so{n}$.\\

Therefore $T$ is a $\so{S}$-storage operator. $\Box$ \\

{\bf Definition} We say that a $\l X$-term $U$ satisfies the property $(P)$ iff for each constant
$X_{l,a,b,\sou{c}}$ of $U$ we have :
\begin{itemize}
\item[] - $a,b,\sou{c}$ satisfy $(P)$
\item[] - $X_{l,a,b,\sou{c}}$ is applied to $a$ and $b$ ;  
\item[] - $a,b$ do not contain free variables which are bounded in $U$.   
\end{itemize}

{\bf Lemma 1} {\it  Let $U,V$ be $\l X$-terms which do not begin by $\l$. If $U$ satisfies $(P)$
and $U \p V$, then $V$ satisfies $(P)$.}\\

{\bf Proof}  It is enough to do the proof for one step of head reduction. We have $U =
(\l x u)v\sou{w}$ and $V = (u[v / x])\sou{w}$. Since $U$ satisfies $(P)$, then $u,v,\sou{w}$
satisfy $(P)$ and $x$ is not free in $a,b$ if the constant $X_{l,a,b,\sou{d}}$ appears in $u$.
Therefore $u[v / x],u_1, ... ,u_m$ satisfy $(P)$ and $V$ satisfies $(P)$. $\Box$ \\

 Let $\D : \L_x \f \L_X$ the simultaneous substitution defined by :
\begin{center}
$\D(x_n)=X_n$\\
for every $0\leq k\leq n-1$, $\D(x_{k,a,b,\sou{c}})=(X_{k,\D(a),\D(b),\sou{\D(c)}})\D(a)\D(b)$\\
$\s(x)=x$ if $x \neq x_n,x_{k,a,b,\sou{c}}$
\end{center}

{\bf Lemma 2} {\it  If $U$ is a $\l X$-term satisfies $(P)$, then there is a $\l x$-term $U'$
such that $\D(U')=U$.}\\

{\bf Proof} By induction on $U$.
\begin{itemize}
\item[] - For $U = x$, it is true.
\item[] - If $U = \l x V$, then $V$ satisfies $(P)$, and, by induction hypothesis, there is a $\l
x$-term $V$ such that $\D(V')=V$. We put $U' = \l x V'$. We have $\D(U')=U$.
\item[] - If $U = (U_1)U_2$ (where $U_1$ does not begin by a constant), then $U_1,U_2$ satisfy
$(P)$, and, by induction hypothesis, there are $\l x$-terms $U'_1,U'_2$ such that $\D(U'_1)=U_1$
and $\D(U'_1)=U_1$. We put $U' = (U'_1)U'_2$. We have $\D(U')=U$. 
\item[] - If $U = (X_{k,a,b,\sou{c}})ab\sou{V}$, then $a,b,\sou{c},\sou{V}$ satisfy
$(P)$, and, by induction hypothesis, there are $\l x$-terms $a',b',\sou{c'},\sou{V'}$ such that
$\D(a')=a$, $\D(b')=b$, $\D(\sou{c'})=\sou{c}$, and $\D(\sou{V'})=\sou{V}$. We put $U' =
(x_{k,a',b',\sou{c'}})\sou{V'}$. We have $\D(U')=U$. $\Box$  
\end{itemize}

{\bf Theorem 2} {\it $T$ is a $\so{S_1}$-storage operator iff $T$ is a storage operator.}\\

{\bf Proof}  Let $n \geq 0$. If $T$ is a $\so{S_1}$-storage operator, then there is a finite
sequence of head reduction $\{ U_i \p V_i \}_{1\leq i\leq r}$ such that :

\begin{itemize}
\item[] 1) $U_i$ and $V_i$ are $\l X$-terms ;
\item[] 2) $U_1 = (T)X_n f$ and $V_r = (f)\t_n$ where $\t_n$ is closed $\l$-term
$\b$-equivalent to $\so{n}$ ;  
\item[] 3) $V_i = (X_n) a b \sou{c}$ or $V_i = (X_{l,a,b,\sou{c}}) u
v \sou{w}$  $0\leq l\leq n-1$ ; 
\item[] 4) If $V_i = (X_n)a b \sou{c}$, then $U_{i+1} = (\so{0})a b \sou{c}$ if $n=0$ and
$U_{i+1} = ((\so{S_1})X_{n-1,a,b,\sou{c}})a b\sou{c}$ if $n \neq 0$ ;
\item[] 5) If $V_i = (X_{l,a,b,\sou{c}})u v\sou{w}$ $0\leq l\leq n-1$, then $U_{i+1} = (\so{0})u v
\sou{w}$ if $l=0$ and $U_{i+1} = ((\so{S_1})X_{l-1,u,v,\sou{w}})u v\sou{w}$ if $l \neq 0$.
\end{itemize}

We prove (by induction on $i$) that, for every $1\leq i\leq r$, $V_i$ satisfies $(P)$.
\begin{itemize} 
\item[] - For $i=1$, it is true. 
\item[] - Assume that is true for $i$, and prove it for $i+1$.
\begin{itemize}
\item[] If $V_i = (X_n)a b \sou{c}$, we have two cases :
\begin{itemize}
\item[] - if $n=0$, then $U_{i+1} = (\so{0})a b \sou{c}$. By induction hypothesis $V_i$ satisfies
$(P)$, then $a, b, \sou{c}$ satisfy $(P)$, therefore $U_{i+1}$ and $V_{i+1}$ satisfy $(P)$.
\item[] - if $n \neq 0$, then $U_{i+1} = ((\so{S})X_{n-1,a,b,\sou{c}})a b\sou{c}$. By induction
hypothesis $V_i$ satisfies $(P)$, then $a, b, \sou{c}$ satisfy $(P)$. Since $U_{i+1} \p
((a)(X_{n-1,a,b,\sou{c}})ab)\sou{c}$, then $V_{i+1}$ satisfies $(P)$. 
\end{itemize}
\item[] If $V_i = (X_{l,a,b,\sou{c}})u v\sou{w}$ $0\leq l\leq n-1$, then $u=a$, $v=b$, and
$\sou{w}=\sou{c}$ since, by induction hypothesis, $V_i$ satisfies $(P)$. We have two cases :
\begin{itemize}
\item[] - if $n=0$, then $U_{i+1} = (\so{0})a b \sou{c}$. By
induction hypothesis $V_i$ satisfies $(P)$, then $a, b, \sou{c}$ satisfy $(P)$, therefore
$U_{i+1}$ and $V_{i+1}$ satisfy $(P)$. 
\item[] - if $n \neq 0$, then $U_{i+1} = ((\so{S})X_{l-1,a,b,\sou{c}})a b\sou{c}$. By induction
hypothesis $V_i$ satisfies $(P)$, then $a, b, \sou{c}$ satisfy $(P)$. Since $U_{i+1} \p
((a)(X_{l-1,a,b,\sou{c}})ab)\sou{c}$, then $V_{i+1}$ satisfies $(P)$.  
\end{itemize}
\end{itemize} 
\end{itemize}

Therefore there is a finite sequence of head reduction $\{ M_i \p N_i \}_{1\leq
i\leq r}$ such that :
\begin{itemize}
\item[] 1) $M_i$ and $N_i$ are $\l X$-terms ;   
\item[] 2) $M_1 = (T)X_n f$ and $N_r = (f)\t_n$ where $\t_n$ is closed $\l$-term
$\b$-equivalent to $\so{n}$ ;  
\item[] 3) $N_i = (X_n) a b \sou{c}$ or $N_i = (X_{l,a,b,\sou{c}}) a b \sou{d}$  $0\leq l\leq n-1$;
\item[] 4) If $N_i = (X_n)a b \sou{c}$, then $M_{i+1} = (b) \sou{c}$ if $n=0$ and
$M_{i+1} = ((a)(X_{n-1,a,b,\sou{c}})a b)\sou{c}$ if $n \neq 0$ ;
\item[] 5) If $N_i = (X_{l,a,b,\sou{c}})a b\sou{d}$ $0\leq l\leq n-1$, then $M_{i+1} =(b)
\sou{d}$ if $l=0$ and $M_{i+1} = ((a)(X_{l-1,a,b,\sou{d}})a b)\sou{d}$ if $l \neq
0$. 
\end{itemize} 

Since, for every $1\leq i\leq r$, $M_i$ and $N_i$ satisfy $(P)$, let $M'_i$ and $N'_i$ the $\l
x$-terms such that : $\D(M'_i)=M_i$ and $\D(N'_i)=N_i$. We have : 
\begin{itemize}
\item[] 1) $M'_i$ and $N'_i$ are $\l x$-terms ;    
\item[] 2) $M'_1 = (T)x_n f$ and $N'_r = (f)\t'_n$ where $\t'_n$ is closed $\l$-term
$\b$-equivalent to $\so{n}$ ;  
\item[] 3) $N'_i = (x_n) a' b' \sou{c'}$ or $N'_i = (x_{l,a',b',\sou{c'}})\sou{d'}$  $0\leq l\leq n-1$;
\item[] 4) If $N'_i = (x_n)a' b' \sou{c'}$, then $M'_{i+1} = (b') \sou{c'}$ if $n=0$ and
$M'_{i+1} = ((a')x_{n-1,a',b',\sou{c'}})\sou{c'}$ if $n \neq 0$ ;
\item[] 5) If $N'_i = (x_{l,a',b',\sou{c'}})\sou{d'}$ $0\leq l\leq n-1$, then $M'_{i+1} =(b')
\sou{d'}$ if $l=0$ and $M'_{i+1} = ((a')x_{l-1,a',b',\sou{d'}})\sou{d'}$ if $l \neq
0$. 
\end{itemize}

Therefore $T$ is a storage operator. $\Box$ \\

{\bf Theorem 3} {\it There is a $\so{S_2}$-storage operator which is a no storage operator.}\\

{\bf Proof} Let $T = \l x (x) ~ a ~ b ~ \so{0} ~ \so{S}$ where \\
$a=\l x \l y \l z ((x)(z)(x)II\l x \so{0})\l x(\so{S})(z)x$, \\
$b=\l x \l y \l z (z)x$, \\
and $\so{S}$ a successor.\\

Let $n$ be an integer.\\
If $n = 0$, then  we check that :
\begin{eqnarray}
(T) ~ X_n ~ f &\p &(X_n) ~ a ~ b ~ \so{0} ~ \so{S} ~ f \nonumber\\
\quad \quad \quad \quad (\so{0}) ~ a ~ b ~ \so{0} ~ \so{S} ~ f &\p &(f)\so{0}\quad \quad \quad \quad \nonumber 
\end{eqnarray}

If $n \neq 0$, then we check that :
\begin{eqnarray}
(T) ~ X_n ~ f &\p &(X_n) ~ a ~ b ~ \so{0} ~ \so{S} ~ f \nonumber\\
((\so{S_2})X_{n-1,a,b,\so{0},\so{S},f}) ~ a ~ b ~ \so{0} ~ \so{S}  ~ f
 &\p &(X_{n-1,a,b,\so{0},\so{S},f}) ~ a ~ (a)b ~ \so{0} ~ \so{S} ~ f \nonumber\\
&. \nonumber\\
&. \nonumber\\
&. \nonumber\\
((\so{S_2})X_{0,a,(a)^{n-1}b,\so{0},\so{S},f}) ~ a ~ (a)^{n-1}b ~ \so{0} ~ \so{S} ~ f
 &\p &(X_{0,a,(a)^{n-1}b,\so{0},\so{S},f}) ~ a ~ (a)^nb ~ \so{0} ~ \so{S} ~ f \nonumber \\
(\so{0}) ~ a ~ (a)^nb ~ \so{0} ~ \so{S} ~ f &\p &(a)^nb ~ \so{0} ~ \so{S} ~ f \nonumber 
\end{eqnarray}

We define two sequences of $\l$-terms $(P_i)_{0 \leq i \leq n}$ and $(Q_i)_{0 \leq i \leq n}$ by
: \begin{center}
$Q_0 = \so{S}$, and, for every $0 \leq k \leq n-1$, we put $Q_{k+1} = \l x (\so{S})(Q_k)x$ \\
$P_0 = \so{0}$, and, for every $0 \leq k \leq n-1$, we put $P_{k+1} = (Q_k)((a)^{n-k-1}b)II\l x \so{0}$
\end{center}

It is easy to check that, for every $1 \leq k \leq n$, $Q_k \simeq_{\b} \l x (\so{S})^{k+1} x$.\\
We prove (by induction on $k$) that, for every $0 \leq k \leq n$, we have $(a)^nb ~ \so{0} ~ \so{S} ~
f \p (a)^{n-k}b ~ P_k ~ Q_k ~ f$. 
\begin{itemize}
\item[] - For $k=0$, it is true. 
\item[] - Assume that is true for $k$, and prove it for $k+1$.\\ 
$(a)^{n-k}b ~ P_k ~ Q_k ~ f = (a) ~ (a)^{n-k-1}b ~ P_k ~ Q_k ~ f \p$ \\
$((a)^{n-k-1}b) ~ (Q_k)((a)^{n-k-1}b)II\l x \so{0} ~ \l x(\so{S})(Q_k)x  ~ f = (a)^{n-{k+1}}b ~ P_{k+1} ~
Q_{k+1} ~ f$.
\end{itemize}
In particular, for $k=n$, we have $(a)^nb ~ \so{0} ~ \so{S} ~ f \p (b) ~ P_n ~ Q_n ~ f \p (f)P_n$.\\
$P_n = (Q_{n-1})(b)II\l x \so{0} \simeq_{\b} (\l x (\so{S})^{n} x)(b)II\l x \so{0} \simeq_{\b}
(\so{S})^{n} (\l x \so{0})I \simeq_{\b} (\so{S})^{n} \so{0} \simeq_{\b} \so{n}$.\\

Therefore $T$ is a $\so{S_2}$-storage operator.\\

We define a sequence of $\l$-terms $(P'_i)_{0 \leq i \leq n}$ by :
\begin{center}
$P'_0 = \so{0}$, and for every $0 \leq k \leq n-1$, we put $P'_{k+1} =
(Q_k)(x_{{n-k-1},a,b,P'_{n-k},Q_{n-k},f})IIJ$ 
\end{center}

We check (as before) that :
\begin{eqnarray}
(T) ~ x_n ~ f &\p &(x_n) ~ a ~ b ~ P'_0 ~ Q_0 ~ f \nonumber \\
(a)x_{n-1,a,b,P'_0,Q_0,f} ~  P'_0 ~ Q_0  ~ f
&\p &(x_{n-1,a,b,P'_0,Q_0,f}) ~ P'_1 ~ Q_1 ~ f \nonumber \\
&. \nonumber \\
&. \nonumber \\
&. \nonumber \\
(a)x_{0,a,b,P'_{n-1},Q_{n-1},f}) ~ P'_{n-1} ~ Q_{n-1} ~ f
&\p &(x_{0,a,b,P'_{n-1},Q_{n-1},f}) ~ P'_n ~ Q_n ~ f \nonumber \\
(b) ~ P'_n ~ Q_n ~ f &\p &(f)P'_n \nonumber 
\end{eqnarray}
But $P'_n = (Q_{n-1})(x_{0,a,b,P'_{n-1},Q_{n-1},f})II\l x \so{0}$ is not closed. \\

Note that $P'_n\simeq_{\b} (\so{S})^{n} (X_{0,a,b,P'_{n-1},Q_{n-1},f})II\l x \so{0} \not \simeq_{\b}
\so{n}$. Indeed, if $(\so{S})^{n} (X_{0,a,b,P'_{n-1},Q_{n-1},f})II\l x \so{0}  \simeq_{\b}
\so{n}$, then $(\so{S})^{n}(\l x_1 \l x_2 \l x_3 (\so{S})\so{0})II\l x \so{0}  \simeq_{\b}
\so{n}$, therefore $\so{n+1} \simeq_{\b} \so{n}$. A contradiction.\\

Therefore $T$ is a no storage operator. $\Box$

\end{document}